\documentclass{amsart}
\usepackage{amssymb}
\usepackage{color}

\newtheorem{theorem}{Theorem}
\newtheorem{lemma}{Lemma}
\newtheorem{proposition}{Proposition}

\newcommand{\pp}{\noindent {\bf Proof. }}
\begin{document}
\title{Codimension growth of solvable Lie superalgebras}

\author[D.D. Repov\v s and M.V. Zaicev]
{Du\v san D. Repov\v s and Mikhail V. Zaicev}

\address{Du\v san D. Repov\v s \\Faculty of Education, and
Faculty of  Mathematics and Physics, University of Ljubljana,
 Ljubljana, 1000, Slovenia}
\email{dusan.repovs@guest.arnes.si}

\address{Mikhail V. Zaicev \\Department of Algebra\\ Faculty of Mathematics and
Mechanics\\  Moscow State University \\ Moscow,119992, Russia}

\email{zaicevmv@mail.ru} 



\begin{abstract}
We study numerical invariants of identities of finite-dimensional solvable Lie superalgebras.
We define new series of finite-dimensional solvable Lie superalgebras $L$ with non-nilpotent derived subalgebra $L'$
and discuss their codimension growth. For the first algebra of this series we prove the
existence and integrality of $exp(L)$.   
\end{abstract}

\date{\today}

\maketitle\vskip 0.2in 

\section{Introduction}    
  
Let $A$ be  an algebra over a field $F$ of characteristic zero. One can define an infinite sequence
$\{c_n(A)\}, n=1,2,\ldots,$ of non-negative integers associated with $A$ called {\em codimension
sequence}. It measures the quantity of polynomial identities of $A$. For many classes of algebras the 
sequence $\{c_n(A)\}$ is exponentially bounded. In particular, this holds for associative PI-algebras
\cite{R}, \cite{L}, for finite-dimensional algebras \cite{BD}, \cite{GZ}, for Kac-Moody Lie algebras
\cite{Z1}, \cite{Z2}, and many others. In this case the sequence $(c_n(A))^{1/n}$ has the lower 
and upper
limits $\underline{exp}(A)$ and the $\overline{exp}(A)$ called the {\em lower} and {\em upper} PI-{\em exponents} of $A$,
respectively. If $\underline{exp}(A)=\overline{exp}(A)$ then there exists an ordinary limit called
the PI-{\em exponent} $exp(A)$ of $A$.
At the end of 1980's Amitsur conjectured that $exp(A)$ exists and is an integer for every associative PI-agebra 
$A$. Amitsur's conjecture was proved in \cite{GZ1}, \cite{GZ2}. Later the existence 
and integrality of PI-exponent was proved for finite-dimensional Lie and Jordan algebras  \cite{GZ}, \cite{GRZ1}, 
\cite{GRZ2}, \cite{GZ4}, \cite{Z3}, \cite{GShZ}. On the other hand, there are 
infinite-dimensional solvable Lie algebras with fractional PI-exponents \cite{ZM}, \cite{VZM}, \cite{M+}.

None of these results can be generalized to Lie superalgebras. There is an infinite series of finite-dimensional 
superalgebras $P(t), t\ge 2,$ where all $P(3), P(4),\ldots$ are simple whereas $P(2)$ is not. For
$L=P(2)$  it was proved in \cite{GZ5} that $exp(L)$ exists and is not an integer. Due to
\cite{GZ5}, there is a serious reason to expect that PI-exponent is fractional for any simple superalgebra
$P(t), t\ge 3$.

For infinite-dimensional Lie superalgebras only some partial results are known \cite{RZPr}, \cite{MZ}. In
particular, in \cite{MZ} it was shown that PI-exponent of a Lie superalgebra $L$ exists and is an
integer, provided that its commutator subalgebra $L^2$ is nilpotent. Note that by the Lie Theorem, the subalgebra
$L^2$ is nilpotent for any finite-dimensional solvable Lie algebra $L$. Unfortunately, finite-dimensional
Lie superalgebras  in general do not satisfy this condition. Hence the result of \cite{MZ} cannot be applied to 
finite-dimensional solvable Lie superalgebras. Although there are examples of finite-dimensional Lie 
superalgebras with the fractional PI-exponent, the following conjecture looks natural: Is it true that any
finite-dimensional solvable Lie superalgebra has an integer exponent?

In this paper we construct new series of finite-dimensional solvable Lie superalgebras $S(t), t=2,3,\ldots$
with non-nilpotent derived subalgebras. For $S(2)$ we prove the existence and integrality of PI-exponent
(Theorem \ref{t1}). We also discuss the following related question concerning graded identities. Every Lie superalgebra $L=L_0\oplus L_1$
is endowed by the natural $\mathbb Z_2$-grading. Hence one can also study asymptotic behavior of graded
codimension sequence $\{c_n^{gr}(L)\}$. It was mentioned in \cite{BD} that $c_n(A)\le c_n^{gr}(A)$ for
any algebra $A=\oplus_{g\in G} A_g$ graded by a finite group $G$. Hence $exp(A)\le exp^{gr}(A)$. In the 
assosciative case there are examples where this inequality is strong. For instance, if $A=F[G]$ is the
group algebra of a finite abelian group $G$ then $exp(A)=1$ whereas $exp^{gr}(A)=|G|$. For Lie superalgebras
similar examples are unknown. On the other hand, there are many examples of simple (associative and nonassociative) algebras with 
$exp^{gr}(A)=exp(A)$. In the present paper we give the first example in the class of solvable Lie superalgebras, 
namely, we prove that $exp^{gr}(S(2))=exp(S(2))$ (Theorem \ref{t2}).

\section{Generalities}

Let $A$ be an algebra over $F$ and let $F\{X\}$ be the absolutely free algebra over $F$ with an
infinite set of generators $X$. A non-associative polynomial $f=f(x_1,\ldots,x_n)\in F\{X\}$ is said
to be an {\em identity} of $A$ if $f(a_1,\ldots,a_n)=0$ for any $a_1,\ldots,a_n\in A$. All identities of $A$ 
form an ideal $Id(A)$ of $F\{X\}$.

Denote by $P_n$ the subspace in $F\{X\}$ of all multilinear polynomials on $x_1,\ldots,x_n\in X$. Then 
$P_n\cap Id(A)$ is the set of all multilinear identities of $A$ of degree $n$. Since ${\rm char}~F=0$,
the sequence of subspaces
$$
\{P_n\cap Id(A)\},~n=1,2,\ldots,
$$
completely defines the ideal $Id(A)$. Denote
$$
P_n(A)=\frac{P_n}{P_n\cap Id(A)}\quad{\rm and}\quad c_n(A)=\dim P_n(A).
$$

The sequence of integers $\{c_n(A)\}, n=1,2,\ldots$, called the {\em codimension sequence} of $A$, is an important
numerical characteristic of $Id(A)$. Analysis of asymptotic behavior of $\{c_n(A)\}$ is one of the
main approaches of the study of identities of algebras.

As it was mentioned in Introduction, there is a wide class of algebras $A$ such that $c_n(A)\le a^n$
for some constant $a$. In this case one can define the {\em lower} and the {\em upper} PI-exponents of $A$ 
as follows:
$$
\underline{exp}(A)=\liminf_{n\to\infty}\sqrt[n]{c_n(A)},\quad
\overline{exp}(A)=\limsup_{n\to\infty}\sqrt[n]{c_n(A)},
$$
respectively. If the ordinary limit exists then we can define the {\em (ordinary)} PI-{\em exponent}
$$
exp(A)=\lim_{n\to\infty}\sqrt[n]{c_n(A)}.
$$

A powerful tool for computing codimensions is the representation theory of symmetric group 
$S_n$. One can define an $S_n$-action on the subspace $P_n$ of multilinear polynomials by setting
$$
\sigma f(x_1,\ldots,x_n)=f(x_{\sigma(1)},\ldots,x_{\sigma(n)})
$$
for $\sigma\in S_n$. Then $P_n$ becomes an $FS_n$-module. Since $P_n\cap Id(A)$ is stable under 
$S_n$-action, then $P_n(A)$ is also an $FS_n$-module and its $S_n$-character 
$$
\chi_n(A)=\chi(P_n(A))
$$
is called the $n$th {\em cocharacter} of $A$. By Mashke's Theorem, $P_n(A)$ is completely reducible, so
\begin{equation}\label{eq1}
\chi_n(A)=\sum_{\lambda\vdash n} m_\lambda\chi_\lambda
\end{equation}
where $\chi_\lambda$ is the irreducible $S_n$-character corresponding to the partition
$\lambda$ of $n$. All details concerning $S_n$-representatiton can be found in \cite{REP}.
The total sum of multiplicities in (\ref{eq1}) is called the $n$th {\em colength} of $A$,
$$
l_n(A)=\sum_{\lambda\vdash n} m_\lambda.
$$
Clearly,
\begin{equation}\label{eq1a}
c_n(A)=\sum_{\lambda\vdash n} m_\lambda d_\lambda
\end{equation}
where $d_\lambda=\deg\chi_\lambda$ is the dimension of the corresponding irreducible representation
and the multiplicities $m_\lambda$ are taken from (\ref{eq1}). It is well-known that the colength
sequence $\{l_n(A)\}$ is polynomially bounded for any finite-dimensional algebra $A$.

\begin{proposition}\label {p1} \cite[Theorem 1]{GMZ}
Let $\dim A=d$. Then
$$
l_n(A)\le d(n+1)^{d^2+d}
$$
for all $n\ge 1$.
\end{proposition}

Throughout the paper we will omit brackets in left-normed products in non-associative algebras, i.e.
$abc=(ab)c, abcd=(abc)d$, etc.

\section{Lie superalgebras $S(t)$}

In this section we introduce an infinite series of finite-dimensional solvable Lie superalgebras
with non-nilpotent commutator subalgebra.

First, let $R$ be an arbitrary associative algebra with involution $*: R\to R$. Consider an associative
algebra $Q$ consisting of $2\times 2$-matrices over $R$
$$
Q = \left\lbrace \left(
           \begin{array}{cc}
             A & B \\
             C & D \\
           \end{array}
         \right)
 \mid A,B,C,D\in R \right\rbrace.
$$
Algebra $Q$ can be naturally endowed by $\mathbb Z_2$-grading $Q=Q_0\oplus Q_1$, where
$$
Q_0 = \left\lbrace \left(
           \begin{array}{cc}
             A & 0 \\
             0 & D \\
           \end{array}
         \right)
\right\rbrace,\quad
Q_1 = \left\lbrace \left(
           \begin{array}{cc}
             0 & B \\
             C & 0 \\
           \end{array}
         \right)
\right\rbrace.
$$
It is well-known that if we define a (super) commutator brackets by setting
$$
[x,y]=xy-(-1)^{|x\|y|}yx
$$
for homogeneous $x,y\in Q_0\cup Q_1$, where $|x|=0$ if $x\in Q_0$ and $|x|=1$ if $x\in Q_1$,
then $Q$ becomes a Lie superalgebra. For basic notions of super Lie theory we refer to 
\cite{Sch}. Denote by
$$
R^+=\{x\in R|~x^*=x\}~,\quad R^-=\{y\in R|~y^*=-y\},
$$
the subspaces of symmetric and skew elements of $R$, respectively. Then the subspace
\begin{equation}\label{eq2}
L = \left\lbrace \left(
           \begin{array}{cc}
             x & y \\
             z & -x^* \\
           \end{array}
         \right)
 \mid x\in R, y\in R^+,z\in R^- \right\rbrace=L_0\oplus L_1
\end{equation}
of $Q$ is a Lie superalgebra under the supercommutator product defined above, where 
even and odd components are
$$
L_0 = \left\lbrace \left(
           \begin{array}{cc}
             x & 0 \\
             0 & -x* \\
           \end{array}
         \right)
\right\rbrace,\quad
L_1 = \left\lbrace \left(
           \begin{array}{cc}
             0 & y \\
             z & 0 \\
           \end{array}
         \right)
\right\rbrace.
$$
Note that if $R=M_t(F)$ is a $t\times t$-matrix algebra, $t\ge 3$, then its subalgebra
$\widetilde L\subset L$ consisting of the matrix
$$
 \left\lbrace \left(
           \begin{array}{cc}
             x & y \\
             z & -x^* \\
           \end{array}
         \right)
\right\rbrace
$$
with traceless matrices $x$ where $x\to x^{*}$  is the transpose involution
is a well-known simple Lie superalgebra $P(t)$ (or $b(t)$ in the notations of \cite{Sch}).

Now we clarify the structure of $R$ in our case. Let $R=UT_t(F)$ be an algebra of 
$t\times t$-upper triangular matrices over $F$. It is well-known (see, for example, \cite{Koshl})
that the reflection across the secondary diagonal is the involution on $R$, hence $L$ 
defined in (\ref{eq2}) is a finite-dimensional Lie superalgebra. We denote this superalgebra
by $S(t)$. Its even component $S_0\simeq UT_t(F)$ is solvable hence the entire $L$ is
also solvable (see, for example, \cite{Sch}). It is not difficult to check that the derived 
subalgebra $L^2$ is not nilpotent and we get the following conclusion.

\begin{proposition}\label{p2}
Let $R$ be the upper triangular $t\times t$-matrix algebra with the involution $*: R\to R$, 
the reflection across the secondary diagonal. Then $S(t)=L=L_0\oplus L_1$ well-defined in (\ref{eq2})
is a finite-dimensional solvable Lie superalgebra, $\dim L=t(t+1)$, with non-nilpotent
commutator subalgebra.
\end{proposition}

Now we will have to deal with the Lie superalgebra $S(2)$. First, we compute supercommutators in the 
associative superalgebra $Q\simeq UT_2(F)\otimes M_2(F)$. If $A,B,C$ and $D$ are 
$2\times 2$-matrices then
\begin{equation}\label{eq3}
\left[ \left(
           \begin{array}{cc}
             A &  0 \\
             0 & -A^* \\
           \end{array}
         \right),
\left(
           \begin{array}{cc}
             0 & B  \\
             0 & 0 \\
           \end{array}
         \right)
 \right] =
 \left(
           \begin{array}{cc}
             0 & AB+BA^* \\
             0 &  0    \\ 
           \end{array}
         \right),
\end{equation}

\begin{equation}\label{eq3a}
\left[ \left(
           \begin{array}{cc}
             A &  0 \\
             0 & -A^* \\
           \end{array}
         \right),
\left(
           \begin{array}{cc}
             0 & 0  \\
             C & 0 \\
           \end{array}
         \right)
 \right] =
 \left(
           \begin{array}{cc}
                 0    &  0 \\
             -A^*C-CA &  0    \\ 
           \end{array}
         \right),
\end{equation}

\begin{equation}\label{eq4}
\left[ \left(
           \begin{array}{cc}
             A &  0 \\
             0 & -A^* \\
           \end{array}
         \right),
\left(
           \begin{array}{cc}
             B &   0  \\
             0 & -B^* \\
           \end{array}
         \right)
 \right] =
 \left(
           \begin{array}{cc}
             AB-BA &     0         \\
               0   & -(AB-BA)^*    \\ 
           \end{array}
         \right),
\end{equation}

\begin{equation}\label{eq4a}
\left[ \left(
           \begin{array}{cc}
             0 &  B \\
             0 &  0 \\
           \end{array}
         \right),
\left(
           \begin{array}{cc}
             0 & 0  \\
             C & 0 \\
           \end{array}
         \right)
 \right] =
 \left(
           \begin{array}{cc}
                 BC &  0 \\
                  0 &  CB \\ 
           \end{array}
         \right).
\end{equation}

From now on, we will not use associative multiplication and will omit square brackets
in the product of elements of Lie superalgebra $S(2)$. That is, $xy=[x,y], xyz=[[x,y],z]$
and so on for $x,y,z\in S(2)$. Let $e_{11}, e_{12}$ and  $e_{22}$ be $2\times 2$-matrix
units. Then $e_{11}^*=e_{22}, e_{22}^*=e_{11}, e_{12}^*=e_{12}$ in $R$ and matrices
$$
a= \left(
           \begin{array}{cc}
               e_{11}-e_{22} &         0      \\
                    0        &   e_{11}-e_{22} \\ 
           \end{array}
         \right),~
b= \left(
           \begin{array}{cc}
                    0 & e_{11}+e_{22} \\
                    0 &       0       \\ 
           \end{array}
         \right),~
c= \left(
           \begin{array}{cc}
                         0        &  0 \\
                    e_{11}-e_{22} &  0 \\ 
           \end{array}
         \right),
$$
$$
d= \left(
           \begin{array}{cc}
               e_{11}+e_{22} &         0       \\
                    0        &  -e_{11}-e_{22} \\ 
           \end{array}
         \right),~
x= \left(
           \begin{array}{cc}
               e_{12} &     0    \\
                  0   &  -e_{12} \\ 
           \end{array}
         \right),~
y= \left(
           \begin{array}{cc}
                    0 & e_{12} \\
                    0 &  0     \\ 
           \end{array}
         \right)
$$
form a basis of $S(2)$. By definition $a,d$ and $x$ are even whereas $b,c$ and $y$ 
are odd.    Using (\ref{eq3}), (\ref{eq3a}), (\ref{eq4}), (\ref{eq4a}) we 
can compute all nonzero products of basis elements,
$$
bc=cb=a,~bd=-db=-2b,~cd=-dc
=2c,~xa=-ax=-2x,
$$
$$
xb=-bx=2y,~ya=-ay=-2y,~yc=cy=-x,
yd=-dy=-2y. 
$$

\section{PI-exponent of $S(2)$}

Since we will have to deal with multialternating sets of
arguments in multilinear and multihomogeneous expressions, it is convenient to use
the following agreement. If $f=f(x_1,\ldots,x_n,y_1,\ldots,y_k)$ is a 
non-associative polynomial, multilinear on $x_1,\ldots,x_n$, then  we denote 
the result of alternation of $f$ on $x_1,\ldots,x_n$ by marking all $x_1,\ldots,x_n$
by the same symbol over $x_i$'s. For example,
$$
\bar x_1y \bar x_2\bar x_3=\sum_{\sigma\in S_3}(sgn~\sigma) x_{\sigma(1)}y  x_{\sigma(2)}  x_{\sigma(3)}
$$
or
$$
(y\bar x_1\widetilde x_1)(\bar x_2\widetilde x_2)(\bar x_3\widetilde x_3)=
$$
$$
\sum_{\sigma\in S_3}\sum_{\tau\in S_3}(sgn~\sigma)(sgn~\tau)(yx_{\sigma(1)}x_{\tau(1)})
(x_{\sigma(2)}x_{\tau(2)})(x_{\sigma(3)}x_{\tau(3)}).
$$

Our next goal is to prove the relation
\begin{equation}\label{eq5}
y(\widetilde b\bar c)(\widetilde c\bar d)(\widetilde d~\bar b)\widetilde a\bar a
=384y.
\end{equation}
Since $aa=ab=ba=ac=ca=ad=da=0$, the left hand side of (\ref{eq5}) is equal to
$y(\widetilde b\bar c)(\widetilde c\bar d)(\widetilde d~\bar b)aa$. Hence it suffices to show
that
\begin{equation}\label{eq6}
y(\widetilde b\bar c)(\widetilde c\bar d)(\widetilde d~\bar b)=96y.
\end{equation}

The left hand side of (\ref{eq6}) can be written as the sum
$$
y(b\bar c)(c\bar d)(d\bar b)+y(c\bar c)(d\bar d)(b\bar b)+y(d\bar c)(b\bar d)(c\bar b)
$$
$$
-y(c\bar c)(b\bar d)(d\bar b)-y(d\bar c)(c\bar d)(b\bar b)-y(b\bar c)(d\bar d)(c\bar b).
$$

Direct computations show that
$$
y(b\bar c)(c\bar d)(d\bar b)=y(bc)(cd)(db)+y(bd)(cb)(dc)=4yacb+4ybac,
$$
$$
y(c\bar c)(d\bar d)(b\bar b)=y(cb)(dc)(bd)+y(cd)(db)(bc)=4yacb+4ycba,
$$
$$
y(d\bar c)(b\bar d)(c\bar b)=y(dc)(bd)(cb)+y(db)(bc)(cd)=4ycba+4ybac,
$$
$$
-y(c\bar c)(b\bar d)(d\bar b)=y(cd)(bc)(db)+y(cb)(bd)(dc)=4ycab+4yabc,
$$
$$
-y(d\bar c)(c\bar d)(b\bar b)=y(db)(cd)(bc)+y(dc)(cb)(bd)=4ybca+4ycab,
$$
$$
-y(b\bar c)(d\bar d)(c\bar b)=y(bd)(dc)(cb)+y(bc)(db)(cd)=4ybca+4yabc.
$$
 Since $yb=0$ and $yab=-2yb=0$, we obtain
$$
y(\widetilde b\bar c)(\widetilde c\bar d)(\widetilde d\bar b)=
8yacb+8ycab+8ycba=-16ycb-8xab-8xba
$$
$$
=16xb+16xb-16ya=96y
$$
and therefore (\ref{eq6}), (\ref{eq5}) hold.

Equality (\ref{eq5})  implies the relation
\begin{equation}\label{eq7}
y\underbrace{(\widetilde b\bar c)(\widetilde c\bar d)(\widetilde d~\bar b)\widetilde a\bar a
\cdots (\widetilde{\widetilde b}\bar{\bar c})(\widetilde{\widetilde c}\bar{\bar d})
(\widetilde{\widetilde d}~\bar{\bar b})\widetilde{\widetilde a}\bar{\bar a}}_m \ne 0
\end{equation}
for any $m\ge 1$. Consider the multilinear polunomial
$$
f_m=w(\widetilde x_1^{(1)} \bar z_1^{(1)}) (\widetilde x_2^{(1)} \bar z_2^{(1)})
(\widetilde x_3^{(1)} \bar z_3^{(1)}) \widetilde x_4^{(1)} \bar z_4^{(1)} 
$$
$$\cdots
(\widetilde{\widetilde x}_1^{(m)} \bar{\bar z}_1^{(m)}) (\widetilde{\widetilde x}_2^{(m)} \bar{\bar z}_2^{(m)})
(\widetilde{\widetilde x}_3^{(m)} \bar{\bar z}_3^{(m)}) \widetilde{\widetilde x}_4^{(m)} \bar{\bar z}_4^{(m)}
$$
of degree $4m+1$. Polynomial $f_m$ depends on $2m$ alternating sets of variables, each of order four.
Moreover, $f_m$ assumes a non-zero value under an evaluation $\varphi: X\to S(2)$ such that
$$
\varphi(w)=y, \varphi(x_1^{(i)})=b, \varphi(x_2^{(i)})=c, \varphi(x_3^{(i)})=d,
\varphi(x_4^{(i)})=a,
$$
$$
\varphi(z_1^{(i)})=c, \varphi(z_2^{(i)})=d, \varphi(z_3^{(i)})=b, \varphi(z_4^{(i)})=a,
i=1,\ldots,m.
$$

Denote $n=8m$ and consider the $S_n$-action on variables $\{x_j^{(i)}, z_j^{(i)}$, $1\le  j \le 4,
1\le i \le m\}$. Under this action the subspace
$$
P_{n+1}=P_{n+1}(w,x_j^{(i)}, z_j^{(i)},  1\le  j \le 4,
1\le i \le m)
$$
becomes an $FS_n$-module. Structure of polynomial $f_m$ and the relation $\varphi(f_m)\ne 0$
show that $e_{T_\lambda}f_m$ is not an identity of $S(2)$, where $e_{T_\lambda}$ is the essential 
idempotent corresponding to some Young tableaux $T_\lambda$  with Young diagram $D_\lambda$
and $\lambda=(2m,2m,2m,2m)\vdash n$. In particular,
\begin{equation}\label{eq8}
c_{n+1}(S(2)) \ge \deg\chi_\lambda.
\end{equation}
From the hook formula for $\deg\chi_\lambda$ and the Stirling formula for factorials it follows
that
\begin{equation}\label{eq9}
\deg\chi_\lambda \ge n^{-5} 4^n,
\end{equation}
provided that $n=8m$ and $\lambda=(2m)^{(4)}$.

Inequalities (\ref{eq8}), (\ref{eq9}) give us the lower bound for codimensions $c_n(S(2))$.

\begin{lemma}\label{l1}
Lower PI-exponent of $S(2)$ satisfies the inequality $\underline{exp}(S(2)) \ge 4$.
\end{lemma}
\pp
Let $n\equiv j(mod~8)$ where $0\le j\le 7$. If $j=1$ then $n=8m+1$ and
$$
c_n(S(2))\ge\frac{4^{n-1}}{(n-1)^5} \ge \frac{1}{5n^5} 4^n
$$
by (\ref{eq8}), (\ref{eq9}). If $j\ne 1$ then there exist $m$ and $1\le i\le 8$ such that
$n=8m+1+i$. In this case the polynomial
$$
g=(e_{T_\lambda}f_m)u_1\cdots u_i
$$
of degree $8m+1+i=n$ is not an identity of $S(2)$ since $\varphi(f_m)=(384)^m y$ for the  
above mentioned evaluation $\varphi$ and $ya=-2y$. Hence
$$
c_n(S(2)) \ge 4^{-8} n^{-5} 4^n.
$$
Therefore $\underline{exp}(S(2)) \ge 4$ and we have completed the proof.
$\Box$

Now we are ready to prove main result of the paper.

\begin{theorem}\label{t1}
PI-exponent of Lie superalgebra $S(2)$ exists  and
$$
exp(S(2))=4.
$$
\end{theorem}
\pp
 The lower bound for $ \underline{exp}(S(2)) \ge 4$ is given by Lemma \ref{l1}.
Hence it suffices to prove the inequality
\begin{equation}\label{eq10}
\overline{exp}(S(2)) \le 4.
\end{equation}
Denote $A=S(2)$ and consider the $n$th cocharacter (\ref{eq1}).

\begin{lemma}\label{l2}
Let $m_\lambda\ne 0$ in (\ref{eq1a}) for $A=S(2)$, $\lambda=(\lambda_1,\ldots,\lambda_k)$. Then
either $k\le 4$ or $k=5$ and $\lambda_5=1$.
\end{lemma}
\pp
Let  $m_\lambda\ne 0$ and $k> 4$. Then there exists Young tableaux $T_\lambda$ such that
$e_{T_\lambda}f\not\in Id(A)$ for some multilinear polynomial $f=f(x_1,\ldots,x_n)$.
Recall that
$$
e_{T_\lambda}=\left(\sum_{\sigma\in R_{T_\lambda}}\sigma\right) 
\left(\sum_{\tau\in C_{T_\lambda}}(sgn~\tau)~\tau\right)
$$
where $R_{T_\lambda}$ and $C_{T_\lambda}$ are the row stabilizer and the column stabilizer of $T_\lambda$
in $S_n$, respectively. Note that the polynomial
$$
g=g(x_1,\ldots,x_n)=\left(\sum_{\tau\in C_{T_\lambda}}(sgn~\tau)~\tau
\right)e_{T_\lambda} f
$$
is also non-identity of $A$. If $k>5$ then $g$ contains an alternating set of variables
$\{x_{i_1},\ldots, x_{i_t}\}$ of order $t\ge 6$. Consider an evaluation $\varphi: X
\to B=\{a,b,c,d,x,y\}$. 
The linear subspace $J=\langle x,y\rangle\subset A$ is a nilpotent ideal of $A$, $J^2=0$. If at
least two of $x_{i_\alpha}$, $1\le\alpha\le  6$, lie in $J$ then $\varphi(g)=0$. But if
$\varphi(x_{i_1}),\ldots,\varphi(x_{i_t})$ take not more than five distinct values in $B$ then
also $\varphi(g)=0$, due to the skew symmetry of $g$. This contradiction shows that $k\le 5$.
Similar arguments imply the restriction $\lambda_5\le 1$ and we have completed the proof of
the lemma.
$\Box$

In light of Lemma \ref{l2}, by Lemma 6.2.4 and Lemma 6.2.5 from \cite{GZbook}, we have
$$
\deg\chi_\lambda < Cn^r 4^n
$$
for some constants $C,r$ if $m_\lambda\ne 0$ in (\ref{eq1a}). Finally, applying Proposition \ref{p1},
 we get the inequality (\ref{eq10}) and the proof is completed.
$\Box$

\section{Graded  PI-exponent of $S(2)$}

Recall the definition of the graded codimension of a $\mathbb Z_2$-graded algebra. Let $A=A_0\oplus A_1$
be an $F$-algebra with $\mathbb Z_2$-grading. Denote by $F\{X,Y\}$ the free algebra on two infinite
sets of generators $X$ and $Y$. Let all $x\in X$ be even and all $y\in Y$ odd. Then this parity
on $X\cup Y$ induces $\mathbb Z_2$-grading on $F\{X,Y\}$. A polynomial $f=f(x_1,\ldots,x_m,
y_1,\ldots,y_n)$ with $x_1,\ldots,x_m\in X,y_1,\ldots,y_n\in Y$ is said to be a {\em graded identity} of 
$A$ if $f=f(a_1,\ldots,a_m, b_1,\ldots,b_n)=0$ for all $a_1,\ldots,a_m\in A_o,b_1,\ldots,b_n\in A_1$.

Given $0\le k\le n$, denote by $P_{k,n-k}$ the subspace of all multilinear polynomials on
$x_1,\ldots,x_k\in X,y_1,\ldots,y_{n-k}\in Y$ and define the integer
$$
c_{k,n-k}(A)=\dim\frac{P_{k,n-k}}{P_{k,n-k}\cap Id^{gr}(A)}
$$
where $Id^{gr}(A)$ is the ideal of graded identities of $A$. Then the value
$$
c_n^{gr}(A)=\sum_{k=0}^n {n\choose k} c_{k,n-k}(A)
$$
is called the {\em graded} $n$th {\em codimension} of $A$. As in the non-graded case, the limits
$$
\underline{exp}^{gr}(A)=\liminf_{n\to\infty}\sqrt[n]{c_n^{gr}(A)},\quad
\overline{exp}^{gr}(A)=\limsup_{n\to\infty}\sqrt[n]{c_n^{gr}(A)},
$$
$$
exp^{gr}(A)=lim_{n\to\infty}\sqrt[n]{c_n^{gr}(A)}
$$
are called the {\em lower}, the {\em upper} and the {\em ordinary graded} PI-{\em exponents} of $A$.

The space $P_{k,n-k}$ has a natural $F[S_k\times S_{n-k}]$-module structure where
the symmetric groups $S_k$ and $S_{n-k}$ act on $\{x_1,\ldots,x_{k}\}$ and on
$\{y_1,\ldots,y_{n-k}\}$, respectively. Since $P_{k,n-k}\cap Ig^{gr}(A)$ is stable under the
$S_k\times S_{n-k}$-action, then the quotient space
$$
P_{k,n-k}(A)=\frac{P_{k,n-k}}{P_{k,n-k}\cap Ig^{gr}(A)}
$$
is also an $F[S_k\times S_{n-k}]$-module and its $S_k\times S_{n-k}$-character has the form
\begin{equation}\label{eq11}
\chi_{k,n-k}(A)=\chi(P_{k,n-k}(A))=\sum_{{\lambda\vdash k\atop \mu\vdash n-k}}
m_{\lambda,\mu}\chi_{\lambda,\mu}.
\end{equation}
In particular,
\begin{equation}\label{eq12}
c_{k,n-k}(A)=\sum_{{\lambda\vdash k\atop \mu\vdash n-k}}
m_{\lambda,\mu}\deg\chi_{\lambda}\deg\chi_\mu.
\end{equation}
The sum of multiplicities
$$
l_n^{gr}(A)=\sum_{k=0}^n\sum_{{\lambda\vdash k\atop \mu\vdash n-k}}
m_{\lambda,\mu}
$$
is called $n$th {\em graded colength} of $A$ and is polynomially bounded if $\dim A<\infty$ (see \cite{Z4})
that is, there are constants $C,r$ such that
\begin{equation}\label{eq13}
l_n^{gr}(A)\le Cn^r.
\end{equation}
Recall that $A_0=\langle a,d,x\rangle , A_1=\langle b,d,y\rangle $ for our superalgebra $A=S(2)$ and $x,y$ belong to nilpotent
ideal $J$, $J^2=0$.

The same argument as in the proof  of Lemma \ref{l2} gives us the following result.

\begin{lemma}\label{l3}
Let $A=S(2)$ and let $m_{\lambda,\mu}\ne 0$ in (\ref{eq11}). Then $\lambda=(\lambda_1)$ or $\lambda=
(\lambda_1,\lambda_2)$ or $\lambda=(\lambda_1,\lambda_2,1)$ and $\mu=(\mu_1)$ or $\mu=(\mu_1,\mu_2)$ or
$\mu=(\mu_1,\mu_2,1)$.
\end{lemma} 

As a consequence of Lemma \ref{l3} and Lemmas 6.2.4, 6.2.5 from \cite{GZbook} we get the following statement.

\begin{lemma}\label{l4}
There are constants $c_,r_0,c_1,r_1$  not depending on $k$ such that
$$
\deg\chi_\lambda \le c_0n^{r_0}2^k,\quad \deg\chi_\mu\le c_1n^{r_1}2^k
$$
for all $\lambda\vdash k,\mu\vdash (n-k)$ if $m_{\lambda,\mu}\ne 0$ in (\ref{eq11}.)
\end{lemma}

Our final result says that $exp(S(2))$ and $exp^{gr}(S(2))$ coincide.

\begin{theorem}\label{t2} $\quad  exp(S(2))=exp^{gr}(S(2))=4$.
\end{theorem}
\pp
It is well-known (see \cite{BD}) that $c_n(A)\le c_n^{gr}(A)$ for any group graded algebra
$A$. Hence, by Theorem \ref{t1},
\begin{equation}\label{eq14}
\underline{exp}^{gr}(S(2))\ge 4.
\end{equation}

Let us prove that
\begin{equation}\label{eq15}
\overline{exp}^{gr}(S(2))\le 4.
\end{equation}
By (\ref{eq13}), Lemma \ref{l3} and Lemma \ref{l4}, we have
$$
c_{k,n-k}(S(2))\le c_3n^{r_3} 2^k 2^{n-k}=c_3n^{r_3} 2^n
$$
for some constants $c_3, r_3$. Then by definition of graded codimensions,
$$
c_n^{gr}(S(2))\le c_3n^{r_3}2^n\sum_k{n\choose k} = c_3n^{r_3} 4^n.
$$
The latter relation proves (\ref{eq15}). Finally, (\ref{eq14})
and
(\ref{eq15})
complete the proof.
$\Box$

\section*{Acknowledgements}
We express our sincere thanks to the referee for the numerous comments and suggestions. 
The first author was supported by the Slovenian Research Agency grants BI-RU/16-18-002,
P1-0292, N1-0083, N1-0064, J1-8131, and J-7025.
The second author was supported by the Russian Science Foundation grant 16-11-10013.

\end{document}